\documentstyle{amsart}

\newtheorem{thm}{Theorem}[section]

\newtheorem{prop}[thm]{Proposition}
\newtheorem{cor}[thm]{Corollary}

\theoremstyle{definition}
\newtheorem{defn}{Definition}

\newtheorem{exmp}{Example}

\theoremstyle{remark}

\newtheorem{rem}{Remark}

\numberwithin{equation}{section}

\newcommand{\C}{\Bbb{C}}

\newcommand{\N}{\Bbb{N}}

\def\Ker{\operatorname {Ker}}
\def\Im{\operatorname {Im}}
\def\End{\operatorname {End}}

\def\Lie{\operatorname {Lie}}

\newcommand{\pair}[2]{\langle {#1},{#2} \rangle}
\newcommand{\hpair}[2]{\langle {#1} \mid {#2} \rangle}

\newcommand{\cgoth}{\frak{c}}
\newcommand{\g}{\frak{g}}
\newcommand{\gl}{\frak{gl}}
\newcommand{\slgoth}{\frak{sl}}
\newcommand{\hgoth}{\frak{h}}
\newcommand{\pgoth}{\frak{p}}
\newcommand{\agoth}{\frak{a}}
\newcommand{\ngoth}{\frak{n}}
\newcommand{\kgoth}{\frak{k}}

\newcommand{\Tb}{\bar{T}}
\newcommand{\Rb}{{\bar{R}}}
\newcommand{\qdet}{\text{det}_q}
\newcommand{\crgl}{\Bbb{C}_R[GL(n)]}
\newcommand{\crsl}{\Bbb{C}_R[SL(n)]}
\newcommand{\cqsl}{\Bbb{C}_q[SL(n)]}
\newcommand{\urgl}{U_R(\frak{gl}(n))}
\newcommand{\ursl}{U_R(\frak{sl}(n))}
\newcommand{\ckR}{\check{R}}
\newcommand{\Uq}{U_q(\g)}

\newcommand{\cqg}{{\Bbb{C}}_q[G]}

\begin{document}


\title{On the Cremmer-Gervais quantizations of $SL(n)$}
\author{Timothy J. Hodges}
\address{University of Cincinnati, Cincinnati, OH 45221-0025, U.S.A.}
\email{hodges@@ucbeh.san.uc.edu}

\maketitle


\begin{abstract}
	The non-standard quantum groups $\crgl$ and $\crsl$ are constructed for a two
parameter version of the Cremmer-Gervais $R$-matrix. An epimorphism is
constructed
from $\crgl$ onto the restricted dual $U_\Rb(\frak{gl}(n-1))$ associated to
a related
smaller $R$-matrix of the same form.

	A related result is proved concerning factorizable Lie bialgebras. For any
such Lie
bialgebra, the dual Lie bialgebra has a canonical homomorphic image which is
again
factorizable.
\end{abstract}


\section{Introduction}

	We define and study some non-standard quantifications $\crsl$ and $\crgl$
of the
algebra of functions on the special and general linear groups. These
algebras are
constructed using a two parameter version of the nonstandard solutions of
the Yang-
Baxter equation discovered by Cremmer and Gervais \cite{CG}. We refer to these
algebras as Cremmer-Gervais quantum groups.
 Our main result states that there is a homomorphism from $\crgl$ onto the
dual, or
quantized universal enveloping algebra, $U_\Rb(\frak{gl}(n-1))$ associated
to a related
smaller $R$-matrix of the same form. This result implies that the algebraic
structure of
$\crsl$ will be significantly different from that of the standard
quantification $\cqsl$.
For instance, in the case when $n=3$, the main result provides a map from
$\Bbb{C}_R[SL(3)]$ onto one of the multi-parameter versions of
$U_q(\frak{gl}(2))$.
Thus while all the simple  finite dimensional modules over
$\Bbb{C}_q[SL(3)]$ are one
dimensional, $\Bbb{C}_R[SL(3)]$ has simple modules of all possible dimensions.

	In contrast to the enormous literature on the standard quantum groups, very
little
attention has been paid to these non-standard examples. In their original
paper, Cremmer
and Gervais define and discuss the bialgebra $A(R)$ associated to their
$R$-matrices.
They also define and discuss the associated symmetric space $S(R)$.
In \cite{BDF1,BDF2}, Balog, Dabrowski and Feh\'{e}r  define the quantum
determinant
and antipode needed for the construction of $\Bbb{C}_R[SL(3)]$ . They also
define and
study the associated quantized universal enveloping algebra
$U_R(\frak{sl}(3))$. The
associated formal deformation was studied by Gerstenhaber, Giaquinto and
Schack in
\cite{GGS}.

	This paper represents another step in the project to understand the algebraic
structure of quantifications of $\C[G]$ for $G$ a connected semi-simple
group. In the
standard and multi-parameter case the classification of the primitive ideals
was achieved in
the series of papers \cite{HL1,HL2,HLT,J1,J2}.
Our main result implies that the primitive spectrum of $\crgl$ contains the
primitive
spectrum of $U_\Rb(\frak{gl}(n-1))$ as a closed subset and hence that any
classification
of primitive ideals of $\crgl$ should include a classification of the
primitive ideals of
$\urgl$.

	In the final section we give a result on the dual of a factorizable Lie
bialgebra that
to some extent `explains' the main results of section 4 and  suggests how
these results
might be broadened to more general quantifications of connected semi-simple
groups
associated to solutions of the modified classical Yang-Baxter equation (MCYBE).
Let $\g$ be a semi-simple complex Lie algebra with a coboundary Lie
bialgebra structure
given by a solution $r$ of the MCYBE and let $\g^*$ be the dual Lie
bialgebra. Then
there is a canonical homomorphic image $\tilde{\g}$ of $\g^*$ which is
reductive and for
which the Lie bialgebra structure is given again by a solution $\tilde{r}$
of the MCYBE.
In particular when $r$ is the `Cremmer-Gervais' classical $r$-matrix for
$\frak{sl}(n)$,
then $\tilde{r}$ is one of the two parameter family of  Cremmer-Gervais
$r$-matrices for
$\frak{gl}(n-1)$.
Thus our main result on Cremmer-Gervais quantum groups may be viewed as a
quantification of this special case of  Theorem \ref{flb}.

	Many of the ideas in this paper were developed in discussions with T.
Levasseur.
We thank him for his very significant contributions.


\section{Braided Hopf algebras and their duals}

	We first view a general context for our main result. We show that the
restricted
dual  $A^\circ$ of a braided Hopf algebra $A$ contains a canonical braided
subalgebra
$U_0$. Hence there is a natural map from $A$ into the dual of $U_0$. In this
and the
next section we consider Hopf algebras over an arbitrary field, $k$.

	  Recall the definition of a braided bialgebra given in \cite{LT}.  (Notice
that this
definition is slightly different from that given in \cite{Do,Ka} and that in
the latter such
algebras are termed {\em cobraided}).

\begin{defn}
	A braiding on a bialgebra $A$ is a bilinear pairing $\hpair{\;}{\;}$ such
that the
following conditions hold for all $a$, $b$, $c$ and $d$ in $A$
\begin{enumerate}
\item $\sum\hpair{a_{(1)}}{b_{(1)}} b_{(2)}a_{(2)} =
\sum\hpair{a_{(2)}}{b_{(2)}}
a_{(1)} b_{(1)}$
\item $\pair{\;}{\;}$ is invertible in $(A \otimes A)^*$;
\item $\hpair{a}{bc} = \sum \hpair{a_{(1)}}{b}\hpair{a_{(2)}}{c}$.
\item $ \hpair{ab}{c} = \sum \hpair{b}{c_{(1)}} \hpair{a}{c_{(2)}}$
\end{enumerate}
\end{defn}

Recall that in this situation, we also have that $\hpair{a}{1}= \epsilon(a)
= \hpair{1}{a}$
for all $a \in A$.

	A braided Hopf algebra is a Hopf algebra which is braided as a bialgebra. The
antipode of a braided Hopf algebra is always bijective by \cite[Theorem
1.3]{Do}.
Moreover, in this case the braiding is $S$-invariant in the sense that
$\hpair{a}{b} =
\hpair{S(a)}{S(b)}$.

	For any braided Hopf algebra $A$, the braiding  induces a pair of Hopf algebra
maps $l^{\pm} : A^{op} \to A^\circ$ given by
$$ 	l^+(a)(b) =  \hpair{a}{b} \quad
\text{and}
\quad 	l^-(a)(b) = \hpair{b}{S(a)} 	$$

Denote by $U^\pm$ the images of $\l^\pm$ respectively. Then $U^\pm$ are
obviously
Hopf subalgebras of $A^\circ$. Define the FRT dual $U(A)$ to be the Hopf
subalgebra
generated by $U^+$ and $U^-$. The braiding on $A$ induces a Hopf pairing on
$U^+
\otimes (U^-)^{op} $ given by $\pi(l^+(a),l^-(b)) = \hpair{a}{b}$ using
which one can
construct the Drinfeld double $U^+ \Join U^-$. The multiplication map
$u\otimes v
\mapsto uv$ is then a Hopf algebra map from $U^+ \Join U^-$ to $U(A)$.
 Set
$$ U_0(A) = U^+ \cap U^-$$
Then $U_0(A)$ is also a Hopf subalgebra of $A^\circ$. When there is no danger
of
ambiguity we shall denote $U_0(A)$ by $U_0$.

\begin{prop} Consider the pairing on $U_0$ induced from the pairing
$\pi^{-1}$ on
$(U^+)^{op} \otimes U^-$. That is, if $x = l^+(a)$ and $y = l^-(b)$ are
elements of
$U^0$, then
$$(x \mid y)  = \hpair{a}{S(b)}$$
This pairing is a braiding on $U_0$.
\end{prop}

\begin{pf}
The last three conditions follow easily from the definition of the pairing.
Condition (1)
follows from the definition of multiplication in the Drinfeld double and the
fact that the
multiplication map is a homomorphism.
\end{pf}

	Thus the FRT dual $U(A)$ of a braided Hopf algebra $A$ contains a canonical
braided Hopf subalgebra $U_0(A)$. Conversely the braided Hopf algebra has a
canonical
homomorphic image which contains the FRT dual of $U_0(A)$:

\begin{thm} \label{Atrans}
	Let $A$ be a braided Hopf algebra and let $B$ be a braided Hopf subalgebra of
$U_0(A)$. There is a natural Hopf algebra map $\phi: A \to B^\circ$.  If $a
\in A$ is such
that $l^\pm(a) \in B$, then
$$\phi(a) = l_B^\pm\circ l_A^\pm(a).$$
Hence $\phi(A) \supset U(B)$.
\end{thm}

\begin{pf} The map $\phi$ is given by $\phi(a)(u) = u(a)$. Suppose that $a
\in A$ is such
that $l^+(a) \in U_0$ and let $u = l^-(b) \in U_0$. Then
$$ \phi(a)(u) = l^-(b)(a) = \hpair{a}{S(b)} = (l^+(a) \mid l^-(b) ) = l_B^+
(l^+_A(a))(u)$$
The result for $l^-(a)$ is proved similarly.
\end{pf}

\begin{exmp}
	In the case of standard quantum groups, these results are not particularly
interesting. Let $G$ be a connected, simply connected semi-simple complex
algebraic
group and let $\g = \Lie G$. Following the notation of \cite{HLT}, let
$\cqg$ and $\Uq$
be the usual quantum group and quantized universal enveloping algebra.
Then the algebra $U_0$ is the usual `Cartan part' $U^0$ with braiding given
by restriction
of the Rosso form. The map $\phi: \cqg \to (U^0)^\circ$ becomes an
identification of the
undeformed algebra of functions on the maximal torus $\C[H]$ with the FRT
dual of
$U^0$.
\end{exmp}

 We shall also need the following presumably well-known observation.

\begin{prop} \label{brmap}
	 Suppose that $A$ and $B$ are braided Hopf algebras. Let $\chi:A \to B$ be a
surjective Hopf algebra map which is braided in the sense that $\pair{a}{a'} =
\hpair{\chi(a)}{\chi(a')}$. The induced map $\chi^*: U(B) \to U(A)$ is an
isomorphism.
\end{prop}


\section{Multi-parameter $R$-matrices and twisted braided bialgebras}

	In this section we outline some results concerning the construction of a
family of
multi-parameter $R$-matrices from a given $R$-matrix. None of the ideas here
are
particularly new. Similar ideas are used in \cite{CR,HLT,Ly,Ma,Re}. The
twists described
below are special cases of the dual of Drinfeld's gauge transformations
\cite{Dr,Ma}. The
cocycle condition ensures that the associator is still trivial.

\begin{defn} Let $A$ be a Hopf algebra. A 2-cocycle on $A$ is an invertible
pairing
$\sigma : A \otimes A \to k$ such that for all $x$, $y$ and $z$ in $A$,
$$\sum \sigma(x_{(1)},y_{(1)}) \sigma(x_{(2)}y_{(2)},z) = \sum
\sigma(y_{(1)},z_{(1)})
	\sigma(x,y_{(2)}z_{(2)})	$$
and $\sigma(1,1) = 1$.
\end{defn}

	Given a 2-cocycle $\sigma$ on a Hopf algebra, one can twist the
multiplication to
get a new Hopf algebra $A_\sigma$. The new multiplication is given by
$$x \cdot y = \sum \sigma(x_{(1)},y_{(1)}) x_{(2)} y_{(2)} \sigma^{-
1}(x_{(3)},y_{(3)}). $$
See \cite{DT} for further details.

\begin{thm}
	Let $A$ be a braided bialgebra with braiding $\beta$. Let $\sigma$ be a
2-cocycle
on $A$. Let $A_\sigma$  be the twisted bialgebra defined above. Then $\sigma
\beta
(\sigma^\top)^{-1}$ (convolution product) is a braiding on $A_\sigma$.
\end{thm}

\begin{pf} One can verify directly that $\sigma \beta (\sigma^\top)^{-1}$
satisfies the
axioms of a braiding. Similar results are proved in  \cite{CR,Ly,Ma,Re}.
\end{pf}

	Suppose that an invertible operator $R: V \otimes V \to V \otimes V$ is a
solution of the Yang Baxter equation
$$R_{12}R_{13}R_{23} = R_{23}R_{13}R_{12}.$$
 Let P be the twist operator $P(v \otimes v') = v'\otimes v$. Then the
operator $RP$ is a
Yang Baxter operator in the sense of \cite{LT}; i.e., an invertible solution of
$$R_{12}R_{23}R_{12} = R_{23}R_{12}R_{23}.$$
We will denote by $A(R)$ the braided bialgebra denoted by $A(RP)$ in \cite{LT}.
Denote by $T_i^k$ the standard generators. The relations in $A(R)$ are then
$$ \sum_{u,v} R_{ji}^{uv}T_u^k T_v^l = \sum_{u,v} T_i^u T_j^v R_{vu}^{kl}$$
and the braiding is given by
$$ \langle T_i^k \mid T_j^l\rangle = R^{lk}_{ji}$$
For a given basis $\{e_i \mid i = 1, \dots, n\}$ for $V$, define $e_{ij}:V
\to V$ by
$e_{ij}(e_k) = \delta_{ik}e_j$. Then we may write $R = \sum R_{ij}^{kl}
e_{ik} \otimes
e_{jl}$.

	Suppose that $\sigma \in (A \otimes A)^*$. Then for any left comodule $V$ over
$A$ there is an induced endomorphism $\sigma_V \in \End V \otimes V$ given by
\begin{equation} \label{sigv}
\sigma_V(v \otimes v') = \sigma(v_{(1)},v'_{(1)}) v_{(2)} \otimes v'_{(2)}
\end{equation}
(see, for example \cite{LT}).

\begin{cor} \label{mpr}
	Let $A=A(R)$ where $R$ is an invertible solution of the Yang Baxter equation.
Let $\sigma$ be a 2-cocycle on $A(R)$. Set $Q=\sigma_V$ and let $R_\sigma =
QR(PQ^{-1}P)$.  Then $R_\sigma$ is a solution of the YBE and $A(R)_\sigma \cong
A(R_\sigma)$ as braided bialgebras.
\end{cor}

\begin{pf} This follows from the universal properties of $A(R)$ given in
\cite{LT}.
\end{pf}

Note that for the corresponding Yang Baxter operators $\ckR_\sigma =
R_\sigma P$ and
$\ckR = RP$ we have the simpler relation $\ckR_\sigma = Q \ckR Q^{-1}$.

\begin{thm}\label{mphr}
	Suppose that $R$ is a solution of the YBE such that
$$R_{ij}^{kl} \neq 0 \Rightarrow i+j=k+l. $$
Let $p\in \C$ and let $Q = \sum p^{j-i}e_{ii} \otimes e_{jj}$. Then $QRQ$ is
also a
solution of the YBE.
\end{thm}

\begin{pf} There is a  bialgebra map from $A(R)$ into the functions on the
2-torus $B=
k[\alpha^{\pm 1},\beta^{\pm 1}]$ given by $T^j_i \mapsto \alpha \beta^i
\delta_{ij}$.
There is a 2-cocycle on $B$ given by $\sigma(\alpha^i \beta^j, \alpha^k
\beta^l)=p^{i+l-j-
k}$ \cite{HLT}. This induces a 2-cocycle on $A(R)$ such that
$\sigma(T_i^k,T_j^l) =
p^{j-i}\delta_{ik}\delta_{jl}$. Since $PQ^{-1}P =Q$, the result follows from
the
previous corollary.
\end{pf}

\begin{rem} Cotta-Ramusino and Rinaldi have a similar  result (Thm 2.1). No
proof is
given and this result appears only to hold for the standard R-matrices.
\end{rem}


\section{Multi-parameter Cremmer-Gervais $R$-matrices}

	We now apply the results of the previous section to construct a two-parameter
family of  solutions of the YBE from the one-parameter family given in
\cite{CG}. For this
and the following section, we shall assume that the base field is the
complex numbers. Let
$0 \neq q,p \in \C$. Set
\begin{align} \label{tprm}
\begin{split}
R=R_n(q,p) =  & \; p^{-1} \left( \sum_{i=1}^{n} q e_{ii} \otimes e_{ii}
 +q\sum_{i>j} p^{-2(i-j)} e_{ii} \otimes e_{jj}
\right.   \\
& + q^{-1}\sum_{i<j} p^{-2(i-j)}  e_{ii} \otimes e_{jj} \\
&   + (q^{-1}-q)\sum_{i<j}\sum_{k=1}^{j-i-1} p^{2k}
e_{i,j-k} \otimes e_{j,i+k}  \\
& + \left.  (q-q^{-1})\sum_{i>j}\sum_{k=0}^{i-j-1} p^{-2k}
e_{i,j+k} \otimes e_{j,i-k} \right)
\end{split}
\end{align}

When $q=p^n$, this is the $R$-matrix given in \cite{CG}. The following
description of
the matrix $R_{ij}^{kl}$ is useful for some of the calculations needed later.
\begin{equation} \label{rfor}
pR_{ij}^{kl} = \begin{cases} qp^{2(l-i)}, & \text{ if } i=k\geq j=l \\
			q^{-1}p^{2(l-i)}, & \text{ if } i=k < j=l\\
			(q-q^{-1})p^{2(l-i)} & \text{ if } j \leq k < i, \; i+j = k+l \\
			(q^{-1}-q)p^{2(l-i)} & \text{ if } i < k < j, \; i+j=k+l
		     \end{cases}
\end{equation}
Notice that $R$ satisfies the homogeneity condition:
$$R_{ij}^{kl} \neq 0 \Rightarrow i+j=k+l. $$
Hence we may apply Theorem \ref{mphr}.

\begin{thm}
	The operator $R=R_n(q,p)$ satisfies the Yang Baxter equation for all $q,p \in
\C^\times$. The Yang Baxter operator $\ckR =pRP$ satisfies the Hecke equation:
$$(\ckR -q)(\ckR +q^{-1})=0.$$
\end{thm}

\begin{pf} Set $Q(p) = \sum p^{j-i} e_{ii} \otimes e_{jj}$. Then for any
$q$, $p$ and
$p'$,
$$p'R(q,pp') = Q(p')R(q,p)Q(p')$$
The fact that $R(q,p)$ satisfies the YBE when $q=p^n$ is proved in
\cite{CG}. It then
follows from \ref{mphr} that $R(q,p)$ satisfies the YBE for all values of
$p$. The fact
that $\ckR$ satisfies the Hecke equation can be calculated directly or
deduced from the
one parameter case using the above identity.
\end{pf}

The corresponding symmetric and exterior algebras are the graded algebras
defined by
$$S(R) = \C\langle x_1,\dots , x_n\rangle / (\{qx_jx_i -pR_{ij}^{kl}x_kx_l\})$$
and
$$\Lambda(R) = \C\langle x_1,\dots , x_n\rangle / (\{q^{-1}x_jx_i
+pR_{ij}^{kl}x_kx_l\})$$
respectively.
Denote by $P_A(t)$ the usual Poincare series of a graded algebra $A$.

\begin{thm} The algebra $\Lambda(R)$ is generated by elements $x_1,\dots,
x_n$ subject
to the relations
$$x_i^2 = 0, \quad x_ix_j =-p^{2(j-i)}x_jx_i$$
The Poincare series of the above algebras are the same as the commutative case:
$$P_S(t) = (1-t)^{-n}, \quad P_\Lambda(t) = (1+t)^n, \quad P_A(t) =
(1-t)^{-n^2}$$
\end{thm}

\begin{pf} In the one parameter case this description of $\Lambda(R)$ is
given in
\cite{HH}. The proof in the multiparameter case is similar. The formula for
$P_\Lambda(t)$ is then clear. The formula for $P_S(t)$ then follows from
\cite{Gu}. The
formula for $P_A(t)$ can be deduced from the one-parameter case \cite{HH} using
\ref{mpr}.
\end{pf}

	This result implies that $A(R)$ has a group-like $q$-determinant element which
can be used to construct the quantum group $\crgl$ and in the one-parameter
case,
$\crsl$. Since $\Lambda^n(R)$ is a one dimensional $A(R)$-comodule generated
by $ x_1
\dots x_n$, there exists a group like element $\qdet$ such that
$$\rho(x_1 \dots x_n) = \qdet \otimes x_1 \dots x_n.$$
where $\rho: \Lambda^n(R) \to A \otimes \Lambda^n(R)$ is the comodule
structure map.

\begin{prop}
\begin{enumerate}
\item $ \qdet = \sum_{\sigma \in S_n} a_\sigma T_1^{\sigma(1)} \dots
T_n^{\sigma(n)}$
where $a_\sigma$ is the scalar such that in $\Lambda(R)$,
$ x_{\sigma(1)} \dots x_{\sigma(n)} = a_\sigma x_1 \dots x_n$.
\item For all $i$ and $k$,
$$\hpair{\qdet}{T_i^k} = (q^{-1}p^n)^{(n-2i)}\delta_{ik}  \quad \text{and}
\quad
\hpair{T_i^k}{\qdet} = (q^{-1}p^n)^{-(n-2i+2)}\delta_{ik}$$
\item For all $i$ and $k$,
$$	T^k_i \qdet = (q^{-1}p^n)^{2(i-k)} \qdet T_i^k 	$$
Hence $\qdet$ is a normal element of $A(R)$.
\end{enumerate}
\end{prop}

\begin{pf}
For the proof of part (1) observe that
\begin{align*}
\rho(x_1 \dots x_n) & = \left(\sum T_1^j \otimes x_j \right) \ldots
\left(\sum T_n^j \otimes
x_j \right) \\
	& = \left( \sum_{\sigma \in S_n} T_1^{\sigma(1)} \ldots T_n^{\sigma(n)}
\otimes
x_{\sigma(1)} \ldots x_{\sigma(n)} \right) \\
		& = \left(\sum _{\sigma \in S_n} a_\sigma T_1^{\sigma(1)} \ldots
T_n^{\sigma(n)} \right) \otimes x_1 \ldots x_n.
\end{align*}
To prove the second assertion in part (2), observe that
\begin{align*}
 \hpair{T_i^k}{\qdet}& = \sum_{\sigma \in S_n} a_\sigma \sum _{j_1,\dots,
j_{n-1}}
	\hpair{T_i^{j_1}}{T_1^{\sigma(1)}} \dots \hpair{T_{j_{n-
1}}^{k}}{T_n^{\sigma(n)}}\\
	&=  \sum_{\sigma \in S_n} a_\sigma \sum _{j_1,\dots, j_{n-1} }
		R_{i1}^{j_1 \sigma(1)} \dots R_{j_{n-1},n}^{k,\sigma(n)}
\end{align*}
Now the homogeneity of $R$ immediately implies that the right hand side is
zero unless $i
= k$.  Assume therefore that $i=k$. Suppose that there exists a $\sigma \neq
e$ and a
multi-index $\{ j_1, \dots j_{n-1} \}$ such that the term
$$		R_{i1}^{j_1 \sigma(1)} \dots R_{j_{n-1},n}^{k,\sigma(n)}	$$
is non-zero. Let $u$ be the smallest integer such that $\sigma(u) > u$. Then
for all $v
<u$, we have $\sigma(v)=v$ and $j_v = i$. Since
$$ \hpair{T^{j_u}_{j_{u+1}}}{T_u^{\sigma(u)}}  \neq 0 	$$
we must have $j_u < i = j_{u-1}$, so $u \leq j_u < i$.  It then follows by
induction that
$j_v \geq u$ for $v \geq u$. Now pick $v >u$ such that $\sigma(v) =u$. Then
$$ 0 \neq \hpair{T^{j_v}_{j_{v-1}}}{T^u_v} \Rightarrow j_{v-1} < u <v $$
contradicting the previous assertion. Thus
$$ \hpair{T_i^k}{\qdet} = \prod_{j=1}^{n} R^{ji}_{ji} = (q^{-1}p^n)^{-(n-
2i+2)}\delta_{ik}$$
as required. The proof of the first assertion is analogous.

	Since $\qdet$ is group-like, the first braiding axiom implies that
$$ \sum \hpair{\qdet}{T_i^m} T^k_m \qdet = \sum \hpair{\qdet}{T_m^k} \qdet
T_i^m
$$
The third assertion then follows from part (2).
\end{pf}

\begin{defn} By the proposition $\qdet$ is a normal element and hence we may
localize
with respect to the Ore set $\{\qdet^k \mid k \in \N \}$.  Denote the algebra
$A(R)[\qdet^{-1}]$ by $\crgl$. In the case where $p^n=q$, $\qdet$ is
central. Denote the
algebra $A(R)/(\qdet-1)$ by $\crsl$.
\end{defn}

\begin{thm}
	The algebra $\crgl$ has a braided Hopf algebra structure such that the
natural map
$A(R) \to \crgl$ is a morphism of braided bialgebras. \par
	Suppose that $p^n=q$. Then $(\qdet -1)$ is a braided Hopf ideal of $\crgl$ and
$\crsl$ is a braided Hopf algebra.
\end{thm}

\begin{pf}
 By \cite[3.1]{Ha}, the localization $A(R)[\qdet^{-1}]$ has a natural
braided bialgebra
structure. On the other hand Gurevitch \cite[5.10]{Gu} proves that this
localization has an
antipode and is therefore a Hopf algebra. The second assertion is clear.
\end{pf}


\section{The main theorem}

	We now investigate the maps defined in the first section in the case where $A=
\crgl$ and $R$ is the 2-parameter Cremmer-Gervais $R$-matrix. In this case
we denote
$U(\crgl)$ by $\urgl$. (Note that since this is the FRT dual, it is not the
standard
definition of a quantized universal enveloping algebra of $\frak{gl}(n)$).
Recall from
\cite{LT} that $\hpair{T_i^k}{T_j^l}= R_{ji}^{lk}$ and
$\hpair{T_i^k}{S(T_j^l)}=
(R^{-1})_{ji}^{lk}$. Hence
$$l^+(T_i^k)(T_j^l) = R_{ji}^{lk}, \qquad l^+(T_i^k)(S(T_j^l)) =
(R^{-1})_{ji}^{lk}$$
and
$$l^-(T_i^k)(T_j^l) = (R^{-1})_{ij}^{kl}, \qquad l^-(T_i^k)(S(T_j^l)) =
R_{ij}^{kl}$$

\begin{thm} \label{rids} The two parameter $R$-matrix $R(q,p)$ satisfies:
\begin{enumerate}
\item $R(q,p)^{-1} = R(q^{-1},p^{-1})$
\item $(R(q,p)^{-1})^{kl}_{ij} = R(q,p)^{l,k+1}_{j,i+1}$
for $i,k <n$ and all $j$ and $l$.
\item $R_n(q,p)_{i+1,j+1}^{k+1,l+1} = R_{n-1}(q,p)_{ij}^{kl}$
\end{enumerate}
\end{thm}

\begin{pf} The Hecke relation
$$\ckR^2 = (q-q^{-1})\ckR + I$$
implies that
$$ R^{-1} = p^2PRP - p(q-q^{-1})P$$
{}From this the first assertion can be proved by direct calculation. The
remaining  assertions
are verified directly using \eqref{rfor}.
\end{pf}

	The following result is given in the case $n=3$ in \cite{BDF1}

\begin{thm} \label{lpm}
\begin{enumerate}
\item For all $i>1$, $l^+(T_i^1) = 0$.
\item For all $i<n$, $l^-(T_i^n) =0$.
\item For all $i,k = 1, \dots, n-1$, $l^-(T_i^k) = l^+(T_{i+1}^{k+1})$.
\end{enumerate}
Thus the elements $t_i^k =l^-(T_i^k) = l^+(T_{i+1}^{k+1})$ belong to $U_0$.
Moreover,
$$(t_i^k\mid t_j^l) = R_{ij}^{kl}$$
for all $i,k = 1, \dots, n-1$.
\end{thm}

\begin{pf} Recall that $\crgl$ is generated by the $T_i^k$ and the $S(T_i^k))$.
For any $j$ and $l$,
$$	l^+(T^1_i)(T^j_l) = R^{l1}_{ji} = 0 \text{ if } i \neq 1$$
and
$$	 l^+(T^1_i)(S(T^j_l)) = (R^{-1})^{l1}_{ji} = 0 \text{ if } i \neq 1$$
An induction argument then shows that  for $i>1$, the $l^+(T^1_i)$  vanish
on $\crgl$,
proving (1).
A similar argument works for part (2) using
$$ 	l^-(T_i^n)(T_j^l) = (R^{-1})^{nl}_{ij} = 0 \text{ if } i \neq n.$$
and
$$ 	l^-(T_i^n)(S(T_j^l)) = (R)^{nl}_{ij} = 0 \text{ if } i \neq n.$$

	To show the equality in part (3), it suffices similarly to show that the
action of the
elements $l^-(T_i^k)$ and $ l^+(T_{i+1}^{k+1})$ coincide on the generators.
Using
Theorem \ref{rids}, we see that
$$	l^+(T^{k+1}_{i+1})(T^l_j) = R^{l,k+1}_{j,i+1}= (R^{-1})^{kl}_{ij} =
l^-(T_i^k)(T_j^l)
$$
and
$$	l^+(T^{k+1}_{i+1})(S(T^l_j)) = (R^{-1})^{l,k+1}_{j,i+1}= R^{kl}_{ij} =
l^-(T_i^k)(S(T_j^l))
$$
as required. Finally,
$$(t_i^k\mid t_j^l) = (l^+(T_{i+1}^{k+1}) \mid l^-(T_j^l ) )
		= \langle T_{i+1}^{k+1} \mid S(T_j^l ) \rangle = (R^{-
1})^{l,k+1}_{j,i+1}
		 = R_{ij}^{kl}.$$
\end{pf}

	Recall that if $R$ is a solution of the Yang-Baxter equation, then so is $R' =
PRP$. Moreover $A(PRP) \cong A(R)^{op}$ as braided bialgebras. Thus if $R$ is a
Cremmer-Gervais $R$-matrix, then we may construct $\C_{R'}[GL(n)]$
analogously and
we will have that $\C_{R'}[GL(n)] \cong \crgl^{op}$.

\begin{thm} \label{into}
Let $R=R_n(q,p)$ and let $\Rb = PR_{n-1}(q,p)P$. There is a homomorphism of
Hopf
algebras $\psi: \C_{\Rb}[GL(n-1)] \to U_R(\gl(n))$ given by
$$	\psi(\Tb_i^k) = l^-(T_i^k) = l^+(T_{i+1}^{k+1})	$$
Moreover, $\Im(\psi) \subset U_0$ and when considered as a map from
$\C_{\Rb}[GL(n-
1)] $ to $U_0$, $\psi$  is a morphism of braided Hopf algebras.
\end{thm}

\begin{pf} Set $t_i^k = l^-(T_i^k) = l^+(T_{i+1}^{k+1})$ for $1\leq i,k \leq
n-1$. Then
$\Delta t_i^k = \sum t_i^j \otimes t_j^k$ and $(t_i^k \mid t_j^l) =
\Rb^{lk}_{ji}$. The
vector space $V=\sum \C t_i^1$ is an $n-1$ dimensional $U_0$-comodule  on
which the
endomorphism \eqref{sigv} induced by the braiding is given by
 $$\sigma_V(t_i^1 \otimes t_j^1) = \bar{R}^{lk}_{ji} \; t_l^1 \otimes t_k^1.$$
 By the universal property of $A(\Rb)$ \cite{LT}, there is a bialgebra map
$\psi: A(\Rb)
\to U_0$ such that
$$	\psi(\Tb_i^k) = t_i^k	$$
Since $\hpair{a}{b} = (\psi(a)\mid \psi(b))$ for all $a$ and $b$ in the
coalgebra spanned
by the $\Tb^k_i$, it follows by induction that $\psi$ is a homomorphism of
braided
bialgebras.
It remains to show that the image of the determinant element of $A(\Rb)$ is
invertible in
$U_0$. For if so, $\psi$ extends to a cobraided bialgebra map $\psi:
\C_{\Rb}[GL(n-1)]
\to U_0$ \cite[3.2]{Ha} which must be a Hopf algebra map \cite[III.8.9]{Ka}.

	Using \ref{lpm} we obtain
$$l^+(\qdet) = \left[ \sum_{\sigma(1)=1} a_\sigma l^+(T^{\sigma(n)}_n) \dots
l^+(T_2^{\sigma(2)}) \right] l^+(T^1_1) $$
and
$$l^-(\qdet) = l^-(T^n_n) \left[  \sum_{\sigma(n)=n} a_\sigma
l^-(T^{\sigma(n-1)}_{n-
1}) \dots l^-(T_1^{\sigma(1)}) \right]  $$
Set $d=  \sum_{\sigma(1)=1} a_\sigma l^+(T^{\sigma(n)}_n) \dots
l^+(T_2^{\sigma(2)})
$ and notice that by \ref{lpm},
$$d =\sum_{\sigma(n)=n} a_\sigma l^-(T^{\sigma(n-1)}_{n-1}) \dots l^-
(T_1^{\sigma(1)}) .$$ Thus $d$ is an invertible element of $U_0$. Denote by
$\overline{\qdet}$ the determinant element of $A(\Rb)$. Then
$$\psi(\overline{\qdet}) =
	\sum_{\sigma \in S_{n-1}}  a_\sigma l^-(T^{\sigma(n-1)}_{n-1}) \dots l^-
(T_1^{\sigma(1)}) =d.$$
\end{pf}

\begin{thm} Let $\phi: \crgl \to \C_{\Rb}[GL(n-1)] ^\circ$ be the induced map.
Then for all $1 \leq i,k \leq n-1$,
$$	\phi(T_{i+1}^{k+1}) = l^+(t^k_i), \qquad \phi(T^k_i) = l^-(t_i^k).	$$
Hence $\phi$ maps onto the subalgebra $U_{\Rb}(\gl(n-1))$.
\end{thm}

\begin{pf} Apply Theorem \ref{Atrans} to the case where $B$ is the image of
$\psi$
above. Since
$l^+(T_{i+1}^{k+1}) = t^k_i$, we obtain that $\phi(T_{i+1}^{k+1}) =
l^+(t^k_i)$.
Similarly one obtains that $\phi(T^k_i) = l^-(t_i^k)$. Since $\phi(T_1^n) =
\phi(T_n^1)
=0$, the image of $\phi$ lies inside $U(B)$. Thus $\phi$ maps onto $U(B)$
which is
isomorphic to $U_{\Rb}(\gl(n-1))$ by Theorem \ref{brmap}.
\end{pf}

	Finally, we interpret the last two results for $\crsl$ in the one parameter
case.
Suppose that $p^n=q$. Set $\ursl = U(\crsl)$ and notice that since $(\qdet
-1)$ is a
braided Hopf ideal, $\ursl$ is naturally isomorphic to $\urgl$.
 Notice also that in this case $\Rb = PR_{n-1}(q,p)P$ is {\em not} the one
parameter
matrix.

\begin{thm}\label{cgsl} Suppose that $q=p^n$. Let $R=R_n(q,p)$ and let $\Rb
= PR_{n-
1}(q,p)P$. There is a homomorphism of braided Hopf algebras $\psi:
\C_{\Rb}[GL(n-1)]
\to \ursl_0$. This map induces a surjective Hopf algebra homomorphism $\phi:
\crsl \to
U_{\Rb}(\gl(n-1))$.
\end{thm}


\section{The dual of a factorizable Lie bialgebra}

	In this final section we prove a result for factorizable (quasi-triangular)
Lie
bialgebras which can be viewed as the `semi-classical limit' of the results
of the previous
section. We prove that, although the dual of a factorizable Lie algebra $\g$
is rarely itself
factorizable, there is a canonical homomorphic image of $\g^*$ that has a
natural
factorizable Lie bialgebra structure. As is often the case, the theorem in
the Lie bialgebra
situation is at the same time much more general and much easier to prove. We
first
formulate the result in such a way as to emphasize the connections with the
results in
Section 2. For the proof and application to reductive Lie algebras, we
reformulate slightly
the notation.

		Recall that a quasi-triangular Lie bialgebra is a Lie bialgebra $\g$ such
that
the cocommutator $\delta$ is of the form $\delta = dr$ where $r\in \g
\otimes \g $ is a
solution of the classical Yang-Baxter equation. Note that the axioms for
$\delta$ imply
that $t= r_{12} + r_{21}$ is an invariant (symmetric) element of $\g \otimes
\g$.
Following \cite{RS} we shall say that $\g$ is factorizable if $t$ is
non-degenerate (as a
bilinear form on $\g^*$). Suppose that $r = \sum r_i \otimes r'_i$. Define
$\phi_\pm : \g^*
\to \g$ by
$$ \phi_+(\xi) = \sum \xi(r_i) r'_i, \quad \phi_-(\xi) = - \sum \xi(r'_i) r_i$$
Then $\phi_\pm$ are Lie algebra homomorphisms whose duals are given by
$\phi_\pm^*
= -\phi_\mp$ \cite{CP,STS}. Recall that for any bialgebra $\g$ we may define
$\g^{op}$
(respectively $\g^{cop}$) to be $\g$ equipped with the same cocommutator but
the
negative or opposite commutator (resp. with the same commutator but the
opposite
cocommutator). Thus we have Lie bialgebra maps,
$$ \phi_\pm : (\g^*)^{cop} \to \g.$$
Let $\cgoth_\pm = \Im \phi_\pm$ and let $\cgoth_0 = \cgoth_+ \cap \cgoth_-$.
Then all
three of $\cgoth_1$, $\cgoth_2$ and $\cgoth_0$ are Lie subbialgebras of
$\g$. Let
$\kgoth= \cgoth^\perp$.
Our main result is the following.

\begin{thm} \label{gflb} The Lie bialgebra $\g^*/\kgoth \cong \cgoth^*$ is
factorizable.
\end{thm}

	Before proving this result, we revert to the notation of \cite{BD}.
Identify $\g^*$
with $\g$ via the bilinear form. The bilinear form on $\g^*$ induces an
invariant bilinear
form on $\g$. Identify $\phi_+$ with  $f\in \End \g$. Then $f$ satisfies
\begin{align}
	 & f+f^* = 1; \label{asy} \\
	 & [f(x),f(y)]= f([x,f(y)]+[f(x),y] -[x,y]). \label{mcy}
\end{align}
(where $f^*$ denote the adjoint of $f$).
Conversely, suppose that $\g$ is a Lie algebra with a fixed  nondegenerate
symmetric
bilinear form. Then  the existence of a factorizable Lie bialgebra structure
on $\g$ is
equivalent to the existence of an operator $f \in \End \g$ such that
\eqref{asy} and
\eqref{mcy} hold.
 We denote the corresponding  factorizable Lie bialgebra by the pair $(\g,f)$.
We may identify $\g^*$ with $\g$ equipped with a new Lie bracket given by
$$[x,y]_f = [x,f(y)]+[f(x),y] -[x,y]$$
(see \cite{STS}).
Denote this Lie algebra by $\g_f$.
As above we have a pair of Lie bialgebra maps,
$$ f, f-1 : \g_f ^{cop} \to \g.$$
Let $\cgoth_1 = \Im( f-1)$ and note that $\cgoth_1^\perp = \Ker f$. So $f$
induces an
endomorphism $\tilde{f}$ of $\cgoth_1/\cgoth_1^\perp$ which satisfies
conditions
\eqref{asy} and \eqref{mcy} above. Hence $\cgoth_1/\cgoth_1^\perp$ has a
factorizable
Lie bialgebra structure.

\begin{thm} The map $f-1: \g_f^{cop} \to \cgoth_1/\cgoth_1^\perp$ is a Lie
bialgebra
homomorphism with kernel $ \kgoth =\Ker f + \Ker (f-1)$.
\end{thm}

\begin{pf} It is clear that $f-1: \g_f^{cop} \to \cgoth_1/\cgoth_1^\perp$ is
a Lie algebra
homomorphism with kernel $ \kgoth =\Ker f + \Ker (f-1)$. Thus it suffices to
show that
the dual map $(f-1)^*:  (\cgoth_1/\cgoth_1^\perp)^* \to (\g_f^{cop})^*$ is
also a Lie
algebra homomorphism. Set $\cgoth = f(\cgoth_1)$. Then $\cgoth^\perp = \kgoth$.
Identifying $(\cgoth_1/\cgoth_1^\perp)^*$ with
$(\cgoth_1/\cgoth_1^\perp)_{\tilde{f}}$
and $(\g_f^{cop})^*$ with $\cgoth^{op}$ via the bilinear form, we find that
$(f-1)^*$
identifies with the vector space isomorphism
$-f: (\cgoth_1/\cgoth_1^\perp)_{\tilde{f}} \to \cgoth^{op}$. Since, for all
$x,y \in
\cgoth_1$,
$$f([\bar{x}, \bar{y}]_{\tilde{f}}) = f([x,y]_f) = [f(x),f(y)]$$
the map $-f$ above is a Lie algebra isomorphism. Hence $f-1$ is an
isomorphism of Lie
bialgebras.
\end{pf}

	Theorm \ref{gflb} then follows immediately from this result using the
appropriate
identifications.

	We now look more closely at the case when $\g$ is a reductive complex Lie
algebra. The work of Belavin and Drinfeld \cite{BD} is easily extended to
give a complete
description of factorizable Lie bialgebra structures. Suppose that $\g$ is a
reductive Lie
algebra equipped with a non-degenerate invariant symmetric bilinear form
$\kappa$. Let
$\hgoth$ be a Cartan subalgebra, let $R$ be the root system and choose a
base $B$ for
$R$. Denote the corresponding positive and negative roots by $R_+$ and $R_-$
respectively. An admissible triple is a triple $(B_1,B_2,\tau)$ where $ B_1,
\; B_2 \subset
B$ and $\tau : B_1 \to B_2$ is a bijection such that
\begin{align}
	& \text{ for all }\alpha, \beta \in  B_1, \;\; (\tau(\alpha),\tau(\beta)) =
(\alpha,
\beta); \\
	& \text{ for all }\alpha \in B_1 \text{ there exists a }t \text{ such that
}\tau^t(\alpha) \not \in B_1.
\end{align}

	For all $\alpha \in R$, choose elements $e_\alpha \in \g^\alpha$ such that
$\kappa(e_\alpha, e_\beta) = \delta_{\alpha,-\beta}$. Set $h_\alpha =
[e_\alpha, e_{-
\alpha}]$. An admissible quadruple is a pair $(\tau, f_0)$ where $f_0 \in
\End \hgoth$
satisfies:
\begin{align}
	& f_0+f_0^* = 1  \\
	& f_0(h_\alpha) = (f_0-1)(h_{\tau(\alpha)}) \text{ for all }\alpha \in B_1
\end{align}
Given any admissible quadruple, we construct an operator $f \in \End (\g)$
in the
following way. First define an ordering on $R_+$ by $\beta \succeq \alpha$
if and only if
there exists a non-negative integer $j$ such that $\beta = \tau^j(\alpha)$.
Define, for all
$\alpha \in R_+$,
$$ f(e_\alpha) = -\sum_{\beta \succ \alpha} e_\beta, \quad
	f(e_{-\alpha}) = \sum_{\beta \preceq \alpha} e_{-\beta}, \quad f|_{\hgoth} =
f_0.$$

	The  proof of Belavin and Drinfeld in the simple case \cite{BD} may be easily
extended to give the following generalization of their original result. (The
result in the
semi-simple case is stated and proved in \cite{CGR}).

\begin{thm}[Belavin-Drinfeld]
	Let $\g$ be a reductive complex Lie algebra equipped with a nondegenerate
invariant symmetric bilinear form. For any admissible quadruple, the
operator $f$ defined
above satisfies \eqref{asy} and \eqref{mcy}. \par
	Conversely, suppose that $f \in \End \g$ satisfies \eqref{asy} and \eqref{mcy}
above.  Then  $f$ is the operator associated to an admissible quadruple for
a suitable
choice of $\hgoth$ and $B$.
\end{thm}

	Given this classification we can describe very precisely the structure of
$(\cgoth_1/\cgoth_1^\perp,\tilde{f})$ in terms of that of $(\g,f)$.

\begin{thm} \label{flb}
 Let $(\tau, f_0)$ be an admissible quadruple and let $(\g,f)$ be the
associated factorizable
Lie bialgebra. Then $\cgoth_1/\cgoth_1^\perp$ is a reductive Lie algebra
with Cartan
subalgebra $\tilde{\hgoth}= \Im (f_0-1)/\Im (f_0-1)^\perp$ and root system
$R_1$. The
operator $\tilde{f}$ is the operator associated to the admissible quadruple
$(\tilde{\tau},
\tilde{f_0})$ where $\tilde{\tau} = \tau|_{\tau^{-1}(B_1\cap B_2)}$ and
$\tilde{f_0}$ is
the map induced by $f_0$ on $\tilde{\hgoth}$.
\end{thm}

\begin{pf} It suffices to recall the description of $\cgoth_1$ and
$\cgoth_1^\perp$ given
in \cite{BD}. Set $\agoth = \sum_{\alpha \in B_1} \C h_\alpha \oplus
\sum_{\alpha \in
R_1} \g^\alpha$; set $\hat{\ngoth} =\sum_{\alpha  \in R_+\backslash R_1}
\g^\alpha$.
Choose $V \subset (\hgoth \cap \agoth)^\perp$ such that $V \supset V^\perp$
and $\Im
(f_0-1) = (\hgoth \cap \agoth) \oplus V$.
Then one verifies easily \cite[\S 6.4]{BD}  that
$$ \cgoth_1 = \pgoth + V, \quad \cgoth_1^\perp = \hat{\ngoth} + V^\perp$$
and
$$\cgoth_1/\cgoth_1^\perp \cong \agoth \oplus V/V^\perp.$$
This provides the first assertion. The second assertion is verified by
calculating $\tilde{f}$
explicitly.
\end{pf}

\begin{exmp} Consider the `Cremmer-Gervais' bialgebra structure on
$\slgoth(n+1)$ (that
is, the bialgebra structure that occurs as the semi-classical limit of the
Cremmer-Gervais
$R$-matrix). Set $\g=\slgoth(n+1)$, $B=\{\alpha_1,\dots ,\alpha_n\}$, $B_1 =
\{\alpha_1, \dots, \alpha_{n-1}\}, B_2=\{\alpha_2, \dots,\alpha_n\}$ and
$\tau(\alpha_i) =
\alpha_{i+1}$. Suppose that $\g$ has the factorizable bialgebra structure
associated to
this triple, (there is a unique choice of $f_0$ in this situation). Then
$\cgoth_1/\cgoth_1^\perp \cong \gl(n)$, and the associated triple is the map
$\tilde{\tau} :
\{\alpha_1, \dots, \alpha_{n-2}\} \to \{\alpha_2, \dots, \alpha_{n-1}\}$
given again by
$\tilde{\tau}(\alpha_i) = \alpha_{i+1}$. Thus we obtain a bialgebra
structure on $\gl(n)$
of `Cremmer-Gervais-type' but with a non-trivial $\tilde{f_0}$. This should
be compared
with Theorem \ref{cgsl}.
\end{exmp}

\end{document}